\documentclass[a4paper,12pt]{article} 
\usepackage{amsmath} 
\usepackage{amsopn} 
\usepackage{amssymb} 
\usepackage[mathscr]{eucal} 
\usepackage{theorem} 
\usepackage{enumerate} 

\theorembodyfont{\itshape} 
\theoremstyle{plain}
  \newtheorem{thm}{Theorem}[section] 
  \newtheorem{pro}[thm]{Proposition} 
   
  \newtheorem{cor}[thm]{Corollary} 

\theorembodyfont{\rmfamily} 
\theoremstyle{plain}
  \newtheorem{defn}{Definition}[section] 
  \newtheorem{ex}{Example}[section] 
  \newtheorem{rem}{Remark}[section] 

\renewcommand{\theequation}% 
           {\thesection.\arabic{equation}}

\begin{document} 

\begin{center} 
{\LARGE Nilpotent structures of oriented neutral vector bundles} 

\vspace{6mm} 

{\Large Naoya {\sc Ando}} 
\end{center} 

\vspace{3mm} 

\begin{quote} 
{\footnotesize \it Abstract} \ 
{\footnotesize In this paper, 
we study nilpotent structures of an oriented vector bundle $E$ 
of rank $4n$ with a neutral metric $h$ and an $h$-connection $\nabla$. 
We define $H$-nilpotent structures of $(E, h, \nabla )$ 
for a Lie subgroup $H$ of $SO(2n, 2n)$ related to 
neutral hyperK\"{a}hler structures. 
We observe that there exist a complex structure $I$ 
and paracomplex structures $J_1$, $J_2$ of $E$ 
such that $h$, $\nabla$, $I$, $J_1$, $J_2$ form 
a neutral hyperK\"{a}hler structure of $E$ 
if and only if there exists an $H$-nilpotent structure of $(E, h, \nabla )$.} 
\end{quote} 

\vspace{3mm} 

\section{Introduction}\label{sect:intro} 

\setcounter{equation}{0} 

Nilpotent structures of an oriented neutral vector bundle $(E, h)$ 
of rank $4n$ are analogues of complex structures and paracomplex structures. 
In particular, 
if $n=1$, 
then a nilpotent structure corresponds to 
a section of one of the light-like twistor spaces associated with $(E, h)$ 
(see \cite{ando7}, \cite{ando8}, \cite{AK} for the light-like twistor spaces) 
and a complex (respectively, paracomplex) structure 
preserving (respectively, reversing) the neutral metric $h$ 
corresponds to a section of 
one of the space-like (respectively, time-like) twistor spaces 
associated with $(E, h)$ (see \cite{ando4}, \cite{BDM2} 
for the space-like twistor spaces and see \cite{ando4}, \cite{HM}, \cite{JR} 
for the time-like twistor spaces). 
A nilpotent structure $N$ gives a null structure on each fiber of $E$ 
and $h$ is null-Hermitian with respect to $N$ (see \cite{dunajski}). 

A nilpotent structure $N$ gives a light-like subbundle $\pi_N$ of $(E, h)$ 
of rank $2n$. 
In addition, $N$ gives a nowhere zero section $\xi_N$ 
of a line bundle $\bigwedge^{2n}\!\pi_N$, 
which is constructed by a local frame field of $\pi_N$ 
given by an admissible frame field of $N$. 
The fact that $\xi_N$ is well-defined is related to a Lie subgroup $G$ 
of $SO(2n, 2n)$. 
The transition function between two admissible frame fields of $N$ is 
valued in $G$. 
It gives a local section of ${\rm End}\,\pi_N$ 
and its determinant is identically equal to one (\cite{ando8}). 
In the present paper, 
we will see that $N$ gives a nondegenerate section $\Theta_N$ 
of $\bigwedge^2\!\pi_N$ satisfying 
$\xi_N =((-1)^{\frac{n(n-1)}{2}} /n!)\Theta^n_N$ 
and that a nondegenerate section $\Theta$ of the $2$-fold exterior power 
of a light-like subbundle $L$ of $(E, h)$ 
of rank $2n$ gives a nilpotent structure $N$ of $(E, h)$ 
satisfying $\pi_N =L$ and $\Theta_N =\Theta$ 
(Theorem~\ref{thm:onetoone}). 
If $n=1$, 
then the section of one of the light-like twistor spaces 
corresponding to $N$ is given by $(1/\sqrt{2} )\xi_N$ (\cite{ando7}). 

Let $\nabla$ be an $h$-connection of $E$, 
i.e., a connection of $E$ satisfying $\nabla h=0$. 
Then the \textit{Walker condition\/} of $N$ is defined, 
based on the definition of Walker manifolds 
(see \cite{BGGNV}, \cite{DDGMMV}, \cite{DDGMMV2}, 
\cite{dunajski}, \cite{walker} for Walker manifolds). 
The Walker condition of $N$ is characterized 
by $\hat{\nabla} \xi_N =\alpha \otimes \xi_N$ 
with a $1$-form $\alpha$ on $M$ 
for the connection $\hat{\nabla}$ of $\bigwedge^{2n}\!E$ 
induced by $\nabla$ ((a) in Proposition~\ref{pro:Wcond}). 
Therefore the Walker condition just means that 
$\hat{\nabla}$ induces a connection of $\bigwedge^{2n}\!\pi_N$. 
If we also denote by $\hat{\nabla}$ the induced connection 
of $\bigwedge^2\!E$, 
then $\nabla N=0$ is equivalent to $\hat{\nabla} \Theta_N =0$, 
and therefore $\nabla N=0$ yields $\hat{\nabla} \xi_N =0$ 
((b) in Proposition~\ref{pro:Wcond}), 
which means that 
a nilpotent structure $N$ parallel with respect to $\nabla$ satisfies 
the Walker condition (see \cite{dunajski}). 
If $n=1$, then $\nabla N=0$ is equivalent to $\hat{\nabla} \xi_N =0$ 
(\cite{ando7}). 
There exists a nilpotent structure $N$ 
satisfying the Walker condition and $\hat{\nabla} \xi_N \not= 0$ 
(Example~\ref{ex:Wnp}, Example~\ref{ex:Wnpn=1}, Example~\ref{ex:Wnpn=1gen}). 
In addition, 
if $n\geq 2$, 
then there exists a nilpotent structure $N$ 
satisfying $\hat{\nabla} \xi_N =0$ and $\nabla N\not= 0$ 
(Example~\ref{ex:dxi0dNnot0}). 
Such an example $N$ 
in Example~\ref{ex:Wnp} $\sim$ Example~\ref{ex:dxi0dNnot0} admits 
a nilpotent structure $N'$ of $(E, h)$ 
satisfying the Walker condition and $E=\pi_N \oplus \pi_{N'}$. 

In the present paper, 
we define special nilpotent structures 
related to a Lie subgroup $K$ of $SO(2n, 2n)$. 
Each of them gives a principal $K$-bundle $P$ associated with $E$ 
and special admissible frame fields we use are given by local sections 
of $P$. 
In addition, $\nabla$ gives a connection in $P$, 
so that the connection form of $\nabla$ with respect to such a local section 
is valued in the Lie algebra of $K$. 
We call such a special nilpotent structure 
a $K$-\textit{nilpotent structure} of $(E, h, \nabla )$. 
If $\nabla$ is flat and if $K\subset G$, 
then for a $K$-nilpotent structure of $(E, h, \nabla )$, 
we can find a local section $e$ of $P$ 
such that each local section of $E$ which appears in $e$ is parallel 
with respect to $\nabla$ (Proposition~\ref{pro:Kflat}). 
A $G$-nilpotent structure of $(E, h, \nabla )$ is just a nilpotent structure 
of $(E, h)$ parallel with respect to $\nabla$ (Proposition~\ref{pro:G2}). 
Let $H$ be a Lie subgroup of $G$ defined as in \eqref{H} below. 
If there exist a complex structure $I$ 
and paracomplex structures $J_1$, $J_2$ of $E$ 
such that $h$, $\nabla$, $I$, $J_1$, $J_2$ form 
a neutral hyperK\"{a}hler structure of $E$, 
then $r(I-(\sin \theta )J_1 +(\cos \theta )J_2 )$ 
($r\in \mbox{\boldmath{$R$}} \setminus \{ 0\}$, $\theta \in [0, 2\pi )$) 
are $H$-nilpotent structures of $(E, h, \nabla )$ (Theorem~\ref{thm:nh}). 
See \cite{BDM2}, \cite{GMV} for paraquaternionic structures, 
and see \cite{DGMY}, \cite{kamada} for neutral hyperK\"{a}hler $4$-manifolds. 
An $H$-nilpotent structure $N$ of $(E, h, \nabla )$ defines 
a unique light-like subbundle $\pi^{\times}_N$ 
of $(E, h)$ of rank $2n$ satisfying $E=\pi_N \oplus \pi^{\times}_N$, 
and a unique $H$-nilpotent structure $N^{\times}$ of $(E, h, \nabla )$ 
(the \textit{dual\/} $H$-\textit{nilpotent structure\/} of $N$) 
satisfying $\pi_{N^{\times}} =\pi^{\times}_N$ (Theorem~\ref{thm:pi'N'}). 
Based on this, we observe that 
if there exists an $H$-nilpotent structure $N$ of $(E, h, \nabla )$, 
then $h$, $\nabla$, 
$I:=(1/2)(N+N^{\times} )$, $J_1 :=-IJ_2$, $J_2 :=(1/2)(N-N^{\times} )$ 
form a neutral hyperK\"{a}hler structure of $E$ (Corollary~\ref{cor:pi'N'}). 

\begin{rem} 
In \cite{ando8}, 
the author studied an $h$-reversing paracomplex structure $J$ of $(E, h)$ 
such that $\nabla J$ is locally represented as $\nabla J=\alpha \otimes N$ 
for a nowhere zero $1$-form $\alpha$ and 
a nilpotent structure $N$ related to $J$. 
The conformal Gauss maps of time-like minimal surfaces in $E^3_1$ give 
such paracomplex structures of the pull-back bundles. 
An oriented neutral $4n$-manifold 
with an almost paracomplex structure as above 
is considered to be a Walker manifold (\cite{ando8}). 
We can find examples of almost paracomplex structures 
of $E^{4n}_{2n}$ as above (\cite{ando8}). 
If $n=1$, 
then an $h$-reversing paracomplex structure of $(E, h)$ as above 
corresponds to a section of one of the time-like twistor spaces 
associated with $(E, h)$ 
such that the covariant derivative is fully light-like (\cite{ando8}). 
We can obtain all the pairs of $h$-reversing almost paracomplex structures 
of $E^4_2$ 
such that each pair gives sections of the two time-like twistor spaces 
with fully light-like covariant derivatives (\cite{ando8}). 
\end{rem} 

\begin{rem} 
The metrics of vector bundles in the present paper are neutral. 
We can refer to \cite{GG} for almost paracomplex structures 
on Riemannian or neutral 4-manifolds. 
See \cite{dunajski0} for a characterization of 
anti-self-dual null-K\"{a}hler $4$-manifolds. 
Neutral metrics appear in the studies of spaces of oriented geodesics 
(\cite{AGK}, \cite{GK}, \cite{salvai}, \cite{salvai2}), 
the ultra-hyperbolic equation 
(\cite{asgeirsson}, \cite{CG}, \cite{guilfoyle}, \cite{john}) 
and quantum field theories (\cite{pavsic}). 
\end{rem} 

This work was supported by 
JSPS KAKENHI Grant Number JP21K03228. 

\section{Nilpotent structures}\label{sect:ans} 

\setcounter{equation}{0} 

Let $M$ be a manifold. 
Let $E$ be an oriented vector bundle over $M$ of rank $4n$ 
and $h$ a neutral metric of $E$. 
Let $N$ be a section of ${\rm End}\,E$. 
We call $N$ a \textit{nilpotent structure\/} of $(E, h)$ 
if on a neighborhood of each point of $M$, 
there exists an ordered pseudo-orthonormal local frame field 
$e=(e_1 , \dots , e_{2n}, e_{2n+1} , \dots , e_{4n} )$ of $(E, h)$ 
satisfying $Ne=e\Lambda_n$ with 
\begin{equation*} 
\Lambda_n :=\left[ \begin{array}{cccc} 
                    O_n & -I_n & O_n &  I_n \\ 
                    I_n &  O_n & I_n &  O_n \\ 
                    O_n &  I_n & O_n & -I_n \\ 
                    I_n &  O_n & I_n &  O_n 
                     \end{array} 
            \right] , 
\end{equation*} 
where $I_n$ is the $n\times n$ unit matrix, 
      $O_n$ is the $n\times n$ zero matrix, 
and we always suppose that $e_1      , \dots , e_{2n}$ are space-like and 
                      that $e_{2n+1} , \dots , e_{4n}$ are time-like. 
Let $N$ be a nilpotent structure of $(E, h)$. 
For $\varepsilon =+$ or $-$, 
we call $N$ an $\varepsilon$-\textit{nilpotent structure\/} of $(E, h)$ 
if on a neighborhood of each point of $M$, 
there exists an ordered pseudo-orthonormal local frame field $e$ 
of $(E, h)$ giving the orientation of $E$ and 
satisfying $NeI'_{4n, \varepsilon} =eI'_{4n, \varepsilon} \Lambda_n$, 
where 
$$I'_{4n, \varepsilon} 
:=\left[ \begin{array}{cccc} 
          I_n & O_n & O_n & O_n \\ 
          O_n & I_n & O_n & O_n \\ 
          O_n & O_n & I_n & O_n \\ 
          O_n & O_n & O_n & I_{n, \varepsilon} 
           \end{array} 
  \right] ,$$ 
and 
\begin{equation*} 
I_{1, \pm} :=\pm 1, \quad 
I_{n, \pm} :=\left[ \begin{array}{cccc} 
                     \pm 1     &   0     & \cdots &  0     \\ 
                         0     &   1     & \ddots & \vdots \\ 
                        \vdots &  \ddots & \ddots &  0     \\ 
                         0     &  \cdots &  0     &  1   
                      \end{array} 
             \right] \ (n\geq 2).  
\end{equation*}
Let $N$ be an $\varepsilon$-nilpotent structure of $(E, h)$. 
Then such a frame field as $e$ is called an \textit{admissible frame field\/} 
of $N$. 
Let $f=(f_1 , \dots , f_{2n} , f_{2n+1} , \dots , f_{4n} )$ be 
an ordered pseudo-orthonormal local frame field of $(E, h)$ 
giving the orientation of $E$. 
Then $f$ is an admissible frame field of $N$ 
if and only if for each admissible frame field $e$ of $N$, 
an $SO(2n, 2n)$-valued function $A$ 
on the intersection of the domains of $e$ and $f$ 
given by $fI'_{4n, \varepsilon} =eI'_{4n, \varepsilon} A$ is valued in 
the Lie subgroup $G$ of $SO(2n, 2n)$ 
defined by $A_0 \Lambda_n =\Lambda_n A_0$ 
for $A_0 \in SO(2n, 2n)$ (\cite{ando8}). 
A section $N$ of ${\rm End}\,E$ is a nilpotent structure of $(E, h)$ 
if and only if $N$ satisfies 
\begin{itemize} 
\item[{\rm (a)}]{${\rm Im}\,N={\rm Ker}\,N$, 
and $\pi_N :={\rm Im}\,N={\rm Ker}\,N$ is a subbundle of $E$ 
of rank $2n$ such that each fiber is light-like,} 
\item[{\rm (b)}]{$h(\phi , N\phi )=0$ for any local section $\phi$ of $E$} 
\end{itemize} 
(\cite{ando7}, \cite{ando8}). 
In particular, $N$ gives a null structure on each fiber of $E$ 
and $h$ is null-Hermitian with respect to $N$ (see \cite{dunajski}). 

Let $E^*$ be the dual bundle of $E$. 
Then $h$ induces a neutral metric $h^*$ of $E^*$. 
For an ordered pseudo-orthonormal local frame field 
$e=(e_1 , \dots , e_{2n}, e_{2n+1} , \dots , e_{4n} )$ of $(E, h)$, 
let $e^* :=(e^1 , \dots , e^{2n}, e^{2n+1} , \dots , e^{4n} )$ be 
the dual frame field of $e$. 
Then $e^*$ is an ordered pseudo-orthonormal local frame field 
of $(E^* , h^* )$. 
Let $N$ be an $\varepsilon$-nilpotent structure of $(E, h)$. 
Then $N$ gives a section $N^*$ of ${\rm End}\,E^*$ by 
\begin{equation} 
N^* \phi^* :=\phi^* \circ N 
\label{N*} 
\end{equation} 
for a local section $\phi^*$ of $E^*$. 
Let $e$ be an admissible frame field of $N$. 
Then $e^*$ satisfies 
$N^* e^* I'_{4n, \varepsilon} =e^* I'_{4n, \varepsilon} {}^t\!\Lambda_n$. 
By this, we obtain $h^* (\phi^* , N^* \phi^* )=0$ 
for any local section $\phi^*$ of $E^*$. 
Therefore, 
if we set 
\begin{equation} 
\Theta_N (\phi^* , \psi^* ):=h^* (\phi^* , N^* \psi^* ) 
\label{Omega0} 
\end{equation}  
for local sections $\phi^*$, $\psi^*$ of $E^*$, 
then $\Theta_N (\psi^* , \phi^* )=-\Theta_N (\phi^* , \psi^* )$. 
This means that $\Theta_N$ is a section of 
the $2$-fold exterior power $\bigwedge^2\!E$ of $E$ given by $N$. 
In addition, 
by $N^* e^* I'_{4n, \varepsilon} =e^* I'_{4n, \varepsilon} {}^t\!\Lambda_n$, 
we obtain 
\begin{equation*} 
\Theta_N =\sum^n_{i=1} \xi_i \wedge \xi_{n+i} , 
\label{Omega} 
\end{equation*} 
where 
\begin{equation} 
  \begin{array}{lcl} 
   \xi_1      :=e_1     -            e_{2n+1} , & \ &  
   \xi_i      :=e_i     -            e_{2n+i} , \\ 
   \xi_{n+1}  :=e_{n+1} +\varepsilon e_{3n+1} , & \ &  
   \xi_{n+i}  :=e_{n+i} +            e_{3n+i}  
    \end{array} \ 
(i=2, \dots , n). 
\label{xie} 
\end{equation} 
Therefore $\Theta_N$ is a section of $\bigwedge^2\!\pi_N$. 
In addition, $\Theta_N$ is nondegenerate, that is, 
$\Theta_N$ gives a nowhere zero section $\xi_N$ of $\bigwedge^{2n}\!\pi_N$ 
by 
\begin{equation} 
\xi_N =\dfrac{(-1)^{\frac{n(n-1)}{2}}}{n!} \Theta^n_N 
      =\xi_1 \wedge \dots \wedge \xi_{2n} . 
\label{xiN}
\end{equation} 

Let $L$ be a subbundle of $E$ of rank $2n$. 
We say that $L$ is an $\varepsilon$-\textit{light-like} subbundle 
of $(E, h)$ if on a neighborhood of each point of $M$, 
there exists an ordered pseudo-orthonormal local frame field 
$e=(e_1 , \dots , e_{2n}, e_{2n+1} , \dots , e_{4n} )$ of $(E, h)$ 
giving the orientation of $E$ 
such that $\xi_1 , \dots , \xi_n , \xi_{n+1} , \dots , \xi_{2n}$ 
as in \eqref{xie} form a local frame field of $L$. 
Let $L$ be an $\varepsilon$-light-like subbundle of $(E, h)$. 
Let $\Theta$ be a nondegenerate section of $\bigwedge^2\!L$. 
Then on a neighborhood of each point of $M$, 
there exists an ordered pseudo-orthonormal local frame field 
$e$ of $(E, h)$ giving the orientation of $E$ and 
satisfying $\Theta =\sum^n_{i=1} \xi_i \wedge \xi_{n+i}$ 
with \eqref{xie}. 
For $\Theta$, let $N^*$ be a section of ${\rm End}\,E^*$ 
given by $\Theta (\phi^* , \psi^* )=h^* (\phi^* , N^* \psi^* )$ 
for local sections $\phi^*$, $\psi^*$ of $E^*$. 
Then $N^*$ satisfies $h^* (\phi^* , N^* \phi^* )=0$ and 
$N^* e^* I'_{4n, \varepsilon} =e^* I'_{4n, \varepsilon} {}^t\!\Lambda_n$. 
For $N^*$, let $N$ be a section of ${\rm End}\,E$ 
given by $N^* \phi^* =\phi^* \circ N$. 
Then $N$ is an $\varepsilon$-nilpotent structure of $(E, h)$ 
such that $e$ is an admissible frame field of $N$ 
and $N$ satisfies $\pi_N =L$ and $\Theta_N =\Theta$. 
Hence we obtain 

\begin{thm}\label{thm:onetoone} 
There exists a one-to-one correspondence 
between the set of $\varepsilon$-nilpotent structures of $(E, h)$ 
and the set of nondegenerate sections of the $2$-fold exterior powers 
of $\varepsilon$-light-like subbundles of $(E, h)$\/$:$ 
each $\varepsilon$-nilpotent structure $N$ corresponds to 
a nondegenerate section $\Theta_N$ of $\bigwedge^2\!\pi_N$ 
by \eqref{N*} and \eqref{Omega0}. 
\end{thm} 

\begin{rem} 
An $\varepsilon$-nilpotent structure $N$ of $(E, h)$ gives 
a nondegenerate section $\Theta_N$ of $\bigwedge^2\!\pi_N$. 
Therefore $\xi_N$ in \eqref{xiN} 
does not depend on the choice of an admissible frame field $e$ of $N$, 
which was already obtained in \cite{ando8}. 
\end{rem} 

\begin{rem} 
Suppose $n=1$. 
The $2$-fold exterior power $\bigwedge^2\!E$ of $E$ is 
a vector bundle over $M$ of rank $6$ 
and $h$ induces a metric $\hat{h}$ of $\bigwedge^2\!E$ of signature (2,4). 
In addition, $\bigwedge^2\!E$ is decomposed as 
$\bigwedge^2\!E=\bigwedge^2_+\!E\oplus \bigwedge^2_-\!E$ 
by two subbundles $\bigwedge^2_+\!E$, $\bigwedge^2_-\!E$ of rank $3$ 
and the restriction of $\hat{h}$ on each of them has signature (1,2). 
The \textit{light-like twistor spaces\/} associated with $(E, h)$ are 
fiber bundles $U_0\!\left(\textstyle\bigwedge^2_{\pm}\!E\right)$ 
in $\bigwedge^2_{\pm}\!E$ respectively such that each fiber is a light cone. 
Each light-like line subbundle 
of $\bigwedge^2_+\!E$ or $\bigwedge^2_-\!E$ 
corresponds to a light-like subbundle of $(E, h)$ of rank $2$ and 
each $\varepsilon$-nilpotent structure $N$ of $(E, h)$ corresponds to 
a section $\Omega_N$ of $U_0\!\left(\bigwedge^2_{\varepsilon}\!E\right)$ 
given by $(1/\sqrt{2} )\xi_N$ (\cite{ando7}, \cite{ando8}). 
\end{rem} 

\begin{rem} 
Suppose $n=1$. 
The space-like twistor spaces 
$U_+\!\left(\textstyle\bigwedge^2_{\pm}\!E\right)$ associated with $(E, h)$ 
are fiber bundles in $\bigwedge^2_{\pm}\!E$ respectively 
such that each fiber is a hyperboloid of two sheets.  
A section of $U_+\!\left(\textstyle\bigwedge^2_{\varepsilon}\!E\right)$ 
corresponds to a complex structure of $E$ preserving $h$. 
See \cite{ando4}, \cite{BDM2} for the space-like twistor spaces. 
The time-like twistor spaces 
$U_-\!\left(\textstyle\bigwedge^2_{\pm}\!E\right)$ associated with $(E, h)$ 
are fiber bundles in $\bigwedge^2_{\pm}\!E$ respectively 
such that each fiber is a hyperboloid of one sheet.  
A section of $U_-\!\left(\textstyle\bigwedge^2_{\varepsilon}\!E\right)$ 
corresponds to a paracomplex structure of $E$ reversing $h$. 
See \cite{ando4}, \cite{HM}, \cite{JR} 
for the time-like twistor spaces. 
See \cite{bryant}, \cite{ES}, \cite{friedrich} for the twistor spaces 
in the case $h$ is a Riemannian (i.e., positive-definite) metric, 
which are the prototypes 
of  $U_+\!\left(\textstyle\bigwedge^2_{\pm}\!E\right)$, 
    $U_-\!\left(\textstyle\bigwedge^2_{\pm}\!E\right)$ 
and $U_0\!\left(\textstyle\bigwedge^2_{\pm}\!E\right)$. 
\end{rem} 

\section{The Walker condition}\label{sect:Wcond} 

\setcounter{equation}{0} 

Let $\nabla$ be an $h$-connection of $E$, i.e., 
a connection of $E$ satisfying $\nabla h=0$. 
Let $N$ be an $\varepsilon$-nilpotent structure of $(E, h)$. 
We say that $N$ satisfies the \textit{Walker condition\/} 
with respect to $\nabla$ 
if for any local section $\psi$ of $\pi_N$, 
$\nabla \psi$ is a $1$-form valued in $\pi_N$. 
See \cite{BGGNV}, \cite{DDGMMV}, \cite{DDGMMV2}, 
\cite{dunajski}, \cite{walker} for Walker manifolds. 

Let $e=(e_1 , \dots , e_{2n}, e_{2n+1} , \dots , e_{4n} )$ be 
an admissible frame field of $N$. 
Let $\omega =[\omega^i_j ]$ be the connection form of $\nabla$ 
with respect to $eI'_{4n, \varepsilon}$. 
Then we have $\nabla eI'_{4n, \varepsilon} =eI'_{4n, \varepsilon} \omega$. 
We represent $\omega$ as 
$$\omega 
 =\left[ \begin{array}{cccc} 
          D_{11} & D_{12} & D_{13} & D_{14} \\ 
          D_{21} & D_{22} & D_{23} & D_{24} \\ 
          D_{31} & D_{32} & D_{33} & D_{34} \\ 
          D_{41} & D_{42} & D_{43} & D_{44} 
           \end{array} 
  \right] ,$$ 
where $D_{kl}$ is the $(k, l)$-block of $\omega$, 
which is an $n\times n$ matrix such that 
the components are given by 
$$\omega^i_j \quad  (i=(k-1)n+1, \dots , kn, \ 
                     j=(l-1)n+1, \dots , ln).$$ 
Since $\omega$ is a $1$-form valued in the Lie algebra of $SO(2n, 2n)$, 
we have 
\begin{equation} 
\begin{split} 
& {}^t\!D_{ii} =-D_{ii} \ (i=1, 2, 3, 4), \quad 
  {}^t\!D_{21} =-D_{12} , \quad 
  {}^t\!D_{43} =-D_{34} , \\ 
& {}^t\!D_{31} = D_{13} , \quad 
  {}^t\!D_{41} = D_{14} , \quad 
  {}^t\!D_{32} = D_{23} , \quad 
  {}^t\!D_{42} = D_{24} . 
\end{split} 
\label{Dij} 
\end{equation} 
In addition, 
$N$ satisfies the Walker condition with respect to $\nabla$ 
if and only if on a neighborhood of each point of $M$, 
an admissible frame field $e$ of $N$ satisfies 
\begin{equation} 
\begin{split} 
D_{11} -D_{13} +D_{31} -D_{33} & =O_n , \\ 
D_{21} -D_{23} -D_{41} +D_{43} & =O_n , \\ 
D_{22} +D_{24} -D_{42} -D_{44} & =O_n .   
\end{split} 
\label{walker} 
\end{equation} 

The connection $\nabla$ induces connections 
of $\bigwedge^2\!E$ and $\bigwedge^{2n}\!E$, 
which are denoted by $\hat{\nabla}$. 
We will prove 

\begin{pro}\label{pro:Wcond} 
Let $N$ be an $\varepsilon$-nilpotent structure of $(E, h)$. 
Then the following hold\/$:$ 
\begin{itemize} 
\item[{\rm (a)}]{$N$ satisfies the Walker condition with respect to $\nabla$ 
if and only if $\hat{\nabla} \xi_N =\alpha \otimes \xi_N$ 
for a $1$-form $\alpha$ on $M;$} 
\item[{\rm (b)}]{$N$ is parallel with respect to $\nabla$ 
if and only if the corresponding section $\Theta_N$ of $\bigwedge^2\!\pi_N$ 
is horizontal with respect to $\hat{\nabla}$, 
and if we suppose one of these conditions, 
then $\hat{\nabla} \xi_N =0$, 
and therefore $N$ satisfies the Walker condition with respect to $\nabla$.} 
\end{itemize} 
\end{pro} 

\vspace{3mm} 

\par\noindent 
\textit{Proof} \ 
If an $\varepsilon$-nilpotent structure $N$ of $(E, h)$ satisfies 
the Walker condition with respect to $\nabla$, 
then we have $\hat{\nabla} \xi_N =\alpha \otimes \xi_N$, 
where $\alpha$ is a $1$-form on $M$ 
which is locally represented as 
\begin{equation} 
\alpha :=\left( -\sum^n_{i=1} \omega^{2n+i}_i 
                +\sum^n_{i=1} \omega^{3n+i}_{n+i}  
         \right) . 
\label{alpha} 
\end{equation} 
If $N$ satisfies $\hat{\nabla} \xi_N =\alpha \otimes \xi_N$ 
for a $1$-form $\alpha$ on $M$, 
then each $\nabla \xi_i$ is a $1$-form valued in $\pi_N$, 
which means that $N$ satisfies the Walker condition. 
Hence we obtain (a) in Proposition~\ref{pro:Wcond}. 
The condition $\nabla N=0$ is equivalent to 
$\omega \Lambda_n =\Lambda_n \omega$. 
This can be rewritten as 
\begin{equation} 
\begin{split} 
& D_{11} -D_{13} = D_{22} -D_{42} = D_{33} -D_{31} = D_{44} -D_{24} , \\ 
& D_{12} +D_{14} =-D_{21} +D_{41} =-D_{23} +D_{43} =-D_{32} -D_{34} . 
\end{split} 
\label{omegaLambda} 
\end{equation} 
We see that \eqref{omegaLambda} is equivalent to $\hat{\nabla} \Theta_N =0$. 
Therefore $\nabla N=0$ is equivalent to $\hat{\nabla} \Theta_N =0$. 
If we suppose $\hat{\nabla} \Theta_N =0$, 
then by \eqref{xiN}, 
we obtain $\hat{\nabla} \xi_N =0$. 
Hence we obtain (b) in Proposition~\ref{pro:Wcond}. 
\hfill 
$\square$ 

\begin{rem} 
It is already known that $\nabla N=0$ means that $N$ satisfies 
the Walker condition with respect to $\nabla$ (see \cite{dunajski}). 
\end{rem} 

\begin{rem} 
Suppose $n=1$. 
Then $\hat{\nabla}$ gives connections of $\bigwedge^2_{\pm}\!E$. 
Let $N$ be an $\varepsilon$-nilpotent structure of $(E, h)$. 
Then the corresponding section 
of $U_0\!\left(\bigwedge^2_{\varepsilon}\!E\right)$ is 
represented as $\Omega_N =(1/\sqrt{2} )\xi_N$, 
and $\nabla N=0$ is equivalent to $\hat{\nabla} \xi_N =0$ (\cite{ando7}). 
Since \eqref{Dij} means 
\begin{equation} 
\omega^i_i     =0 \ (i=1, 2, 3, 4), \quad 
\omega^{i+1}_i =-\omega^i_{i+1} \ (i=1, 3), \quad 
\omega^i_j     =\omega^j_i \ (|i-j|>1), 
\label{Dijn=1} 
\end{equation} 
the Walker condition is given by 
$\omega^1_2 +\omega^1_4 +\omega^3_2 +\omega^3_4 =0$ 
for an admissible frame field $e$ of $N$. 
In addition, noticing \eqref{alpha}, 
we see that $N$ satisfies $\nabla N=0$ if and only if 
not only $\omega^1_2 +\omega^1_4 +\omega^3_2 +\omega^3_4 =0$ 
but also $\omega^3_1 =\omega^4_2$ hold. 
\end{rem} 

\begin{rem} 
Suppose $n=1$. 
For a complex structure $I$ of $E$ preserving $h$, 
$\nabla I=0$ is equivalent to $\hat{\nabla} \Omega_I =0$ 
for the section $\Omega_I$ 
of $U_+\!\left(\textstyle\bigwedge^2_+\!E\right)$ 
or $U_+\!\left(\textstyle\bigwedge^2_-\!E\right)$ 
corresponding to $I$ (\cite{ando4}). 
Similarly, for a paracomplex structure $J$ of $E$ reversing $h$, 
$\nabla J=0$ is equivalent to $\hat{\nabla} \Omega_J =0$ 
for the section $\Omega_J$ 
of $U_-\!\left(\textstyle\bigwedge^2_+\!E\right)$ 
or $U_-\!\left(\textstyle\bigwedge^2_-\!E\right)$ 
corresponding to $J$ (\cite{ando4}). 
\end{rem} 

\begin{rem}\label{rem:J} 
Let $E$, $h$ be as in the beginning of Section~\ref{sect:ans}. 
Let $J$ be a section of ${\rm End}\,E$. 
As in \cite{ando8}, 
we say that $J$ is an $\varepsilon$-\textit{paracomplex structure\/} 
of $(E, h)$ if $J$ satisfies 
\begin{itemize} 
\item[{\rm (i)}]{$J$ is a paracomplex structure of $E$,}
\item[{\rm (ii)}]{$J$ is $h$-reversing, that is, $J^* h=-h$,} 
\item[{\rm (iii)}]{on a neighborhood of each point of $M$, 
there exists an ordered pseudo-orthonormal local frame field 
$e=(e_1 ,\dots e_{2n} , e_{2n+1} , \dots , e_{4n} )$ of $E$ 
giving the orientation and satisfying 
\begin{equation} 
JeI'_{4n, \varepsilon} 
=eI'_{4n, \varepsilon} \Lambda_{n, -} , \quad 
  \Lambda_{n, -} 
:=\left[ \begin{array}{cccc} 
          O_n &  O_n & I_n &  O_n \\ 
          O_n &  O_n & O_n & -I_n \\ 
          I_n &  O_n & O_n &  O_n \\ 
          O_n & -I_n & O_n &  O_n 
           \end{array} 
  \right] . 
\label{J} 
\end{equation}} 
\end{itemize}  
Let $J$ be an $\varepsilon$-paracomplex structure of $(E, h)$. 
Then such a frame field as $e$ is called 
an \textit{admissible frame field\/} of $J$. 
We see that $J$ gives a section $J^*$ of ${\rm End}\,E^*$ 
by $J^* \phi^* =\phi^* \circ J$ 
for a local section $\phi^*$ of $E^*$. 
Then we have 
$J^* e^* I'_{4n, \varepsilon} =e^* I'_{4n, \varepsilon} \Lambda_{n, -}$. 
Therefore, 
if we set $\Theta_J (\phi^* , \psi^* ):=h^* (\phi^* , J^* \psi^* )$, 
then $\Theta_J$ is a section of $\bigwedge^2\!E$ and 
locally represented as 
\begin{equation} 
 \Theta_J 
=\sum^n_{i=1} (e_i      \wedge e_{2n+i} 
                      +(\varepsilon 1)^{\delta_{i1}} 
               e_{3n+i} \wedge e_{n+i} ). 
\label{ThetaJ} 
\end{equation} 
Let $\nabla$ be an $h$-connection of $E$. 
Then $\nabla J=0$ is equivalent to $\hat{\nabla} \Theta_J =0$. 
Suppose that $\nabla J$ is locally represented as the tensor product 
of a $1$-form $\alpha$ and an $\varepsilon$-nilpotent structure $N$ 
so that $e$ is an admissible frame field of both $J$ and $N$. 
Then 
\begin{equation} 
\begin{split} 
&                                \omega^{n+i}_j      
  +(\varepsilon 1)^{\delta_{i1}} \omega^{3n+i}_{2n+j} 
  = \omega^{n+i}_{2n+j} 
  +(\varepsilon 1)^{\delta_{i1}} \omega^{3n+i}_j  
  = \alpha \delta^i_j , \\ 
&   \omega^i_j          -\omega^{2n+i}_{2n+j} 
  = \omega^i_{2n+j}     -\omega^{2n+i}_j 
  = 0, \\ 
&   \omega^{n+i}_{n+j}  
  -(\varepsilon 1)^{\delta_{i1} +\delta_{j1}} \omega^{3n+i}_{3n+j} 
  = \omega^{n+i}_{3n+j} 
  -(\varepsilon 1)^{\delta_{i1} +\delta_{j1}} \omega^{3n+i}_{n+j}  
  = 0     
\end{split} 
\label{omega}
\end{equation}
are obtained in \cite{ando8} for $i, j=1, \dots , n$. 
Therefore by \eqref{ThetaJ} and \eqref{omega}, 
we obtain $\hat{\nabla} \Theta_J =\alpha \otimes \Theta_N$. 
\end{rem} 

Let $N$, $N'$ be $\varepsilon$-nilpotent structures of $(E, h)$. 
Suppose that on a neighborhood of each point of $M$, 
there exists an ordered pseudo-orthonormal local frame field 
$e=(e_1 , \dots , e_{2n}, e_{2n+1} , \dots , e_{4n} )$ of $(E, h)$ 
giving the orientation of $E$ and satisfying 
\begin{itemize} 
\item[{\rm (i)}]{$e$   is an admissible frame field of $N$,} 
\item[{\rm (ii)}]{$e':=(e_1 , \dots , e_{2n}, -e_{2n+1} , \dots , -e_{4n} )$ 
is an admissible frame field of $N'$.} 
\end{itemize} 
Then $E=\pi_N \oplus \pi_{N'}$. 
Referring to \eqref{walker}, we obtain 

\begin{pro}\label{pro:bW} 
Both of $N$ and $N'$ satisfy the Walker condition 
if and only if 
the connection form of $\nabla$ with respect to $e$ as above satisfies 
not only \eqref{Dij} but also 
\begin{equation} 
\begin{split} 
& D_{11} = D_{33} , \quad 
  D_{22} = D_{44} , \quad 
  D_{13} = D_{31} , \quad 
  D_{24} = D_{42} , \\ 
& D_{43} =-D_{21} , \quad 
  D_{41} =-D_{23} . 
\end{split} 
\label{NN'W} 
\end{equation} 
\end{pro} 

\begin{ex}\label{ex:Wnp} 
Suppose that $M$ is diffeomorphic to $\mbox{\boldmath{$R$}}^m$ 
($m\geq 2$) and that an $h$-connection $\nabla$ of $E$ is flat. 
Then the curvature tensor of $\nabla$ vanishes. 
Let $a_1$, $b_1$, $a_2$, $b_2$ be functions on $M$. 
Suppose $n\geq 2$ and 
referring to \cite{ando8}, 
let $F_k$ ($k=1, 2$) be $n\times n$ symmetric matrices defined by 
\begin{equation*} 
F_k :=\left[ \begin{array}{ccccc} 
              a_k   &  b_k   &   0    & \cdots &  0     \\ 
              b_k   &  a_k   &  b_k   & \ddots & \vdots \\ 
              0     &  b_k   & \ddots & \ddots &  0     \\ 
             \vdots & \ddots & \ddots & \ddots &  b_k   \\ 
              0     & \cdots &   0    &  b_k   &  a_k 
               \end{array} 
      \right] . 
\end{equation*} 
Then we have $dF_k \wedge dF_k =O_n$ for $k=1, 2$. 
We set 
\begin{equation*} 
 \omega 
=\left[ \begin{array}{cccc} 
          O_n &  O_n & dF_1 &  O_n \\ 
          O_n &  O_n &  O_n & dF_2 \\ 
         dF_1 &  O_n &  O_n &  O_n \\ 
          O_n & dF_2 &  O_n &  O_n 
          \end{array} 
 \right] . 
\end{equation*} 
Then we have $d\omega +\omega \wedge \omega =O_{4n}$. 
Since we suppose that $\nabla$ is flat, 
there exists an ordered pseudo-orthonormal frame field 
$e=(e_1 , \dots , e_{2n}, e_{2n+1} , \dots , e_{4n} )$ of $(E, h)$ 
giving the orientation of $E$ and 
satisfying $\nabla eI'_{4n, \varepsilon} =eI'_{4n, \varepsilon} \omega$. 
Then $\omega$ is the connection form of $\nabla$ 
with respect to $eI'_{4n, \varepsilon}$ 
and it satisfies \eqref{Dij} and \eqref{NN'W}. 
Let $N$ be an $\varepsilon$-nilpotent structure of $(E, h)$ 
such that $e$ is an admissible frame field of $N$, 
and $N'$   an $\varepsilon$-nilpotent structure of $(E, h)$ 
such that $e'$ is an admissible frame field of $N'$. 
Then $N$, $N'$ satisfy the Walker condition. 
We have 
\begin{equation*} 
\hat{\nabla} \xi_N    =n(-da_1 +da_2 )\otimes \xi_N , \quad 
\hat{\nabla} \xi_{N'} =n( da_1 -da_2 )\otimes \xi_{N'} . 
\end{equation*} 
Therefore if $da_1 \not= da_2$, 
then $\hat{\nabla} \xi_N$, $\hat{\nabla} \xi_{N'}$ do not become zero 
and therefore by (b) in Proposition~\ref{pro:Wcond}, 
$N$, $N'$ are not parallel with respect to $\nabla$. 
\end{ex} 

\begin{ex}\label{ex:Wnpn=1} 
Let $M$, $\nabla$ be as in Example~\ref{ex:Wnp}. 
Let $a_1$, $a_2$ be functions on $M$ satisfying $da_1 \not= da_2$. 
Suppose $n=1$ and set 
\begin{equation*} 
 \omega 
=\left[ \begin{array}{cccc} 
          0   &  0   & da_1 &  0   \\ 
          0   &  0   &  0   & da_2 \\ 
         da_1 &  0   &  0   &  0   \\ 
          0   & da_2 &  0   &  0 
          \end{array} 
 \right] . 
\end{equation*} 
Then we have $d\omega +\omega \wedge \omega =O_4$. 
Therefore there exists an ordered pseudo-orthonormal frame field 
$e=(e_1 , e_2 , e_3 , e_4 )$ of $(E, h)$ 
giving the orientation of $E$ and 
satisfying $\nabla eI'_{4, \varepsilon} =eI'_{4, \varepsilon} \omega$. 
Then $\omega$ is the connection form of $\nabla$ 
with respect to $eI'_{4, \varepsilon}$ 
and it satisfies \eqref{Dij} and \eqref{NN'W}. 
Let $N$, $N'$ be as in Example~\ref{ex:Wnp}. 
Then $N$, $N'$ satisfy the Walker condition. 
We have 
\begin{equation*} 
\hat{\nabla} \xi_N    =(-da_1 +da_2 )\otimes \xi_N , \quad 
\hat{\nabla} \xi_{N'} =( da_1 -da_2 )\otimes \xi_{N'} . 
\end{equation*} 
Therefore $\hat{\nabla} \xi_N$, $\hat{\nabla} \xi_{N'}$ do not become zero 
and therefore $N$, $N'$ are not parallel with respect to $\nabla$. 
\end{ex} 

\begin{ex}\label{ex:Wnpn=1gen}  
Let $M$, $\nabla$ be as in Example~\ref{ex:Wnp} 
and suppose $n=1$. 
Let $N$ be an $\varepsilon$-nilpotent structure of $E$ 
with the Walker condition 
and $e$ an admissible frame field of $N$ defined on $M$. 
Then the connection form $\omega =[\omega^i_j ]$ of $\nabla$ 
with respect to $eI'_{4, \varepsilon}$ 
satisfies $d\omega +\omega \wedge \omega =0$ and 
$\omega^1_2 +\omega^1_4 +\omega^3_2 +\omega^3_4 =0$. 
We can rewrite $d\omega +\omega \wedge \omega =0$ into 
\begin{equation} 
\begin{split} 
d\omega^2_1 +\omega^3_2 \wedge \omega^3_1 +\omega^4_2 \wedge \omega^4_1 & =0, 
 \\ 
d\omega^3_1 +\omega^3_2 \wedge \omega^2_1 -\omega^4_3 \wedge \omega^4_1 & =0, 
 \\ 
d\omega^4_1 +\omega^4_2 \wedge \omega^2_1 +\omega^4_3 \wedge \omega^3_1 & =0, 
 \\ 
d\omega^3_2 -\omega^3_1 \wedge \omega^2_1 -\omega^4_3 \wedge \omega^4_2 & =0, 
 \\ 
d\omega^4_2 -\omega^4_1 \wedge \omega^2_1 +\omega^4_3 \wedge \omega^3_2 & =0, 
 \\ 
d\omega^4_3 +\omega^4_1 \wedge \omega^3_1 +\omega^4_2 \wedge \omega^3_2 & =0. 
\label{domegaij} 
\end{split} 
\end{equation} 
Using the second and fifth equations in \eqref{domegaij} 
and $\omega^1_2 +\omega^1_4 +\omega^3_2 +\omega^3_4 =0$, 
we find a function $a$ satisfying 
\begin{equation} 
da =\omega^3_1 -\omega^4_2 ; 
\label{dalpha} 
\end{equation} 
using the first and sixth equations in \eqref{domegaij} 
and $\omega^1_2 +\omega^1_4 +\omega^3_2 +\omega^3_4 =0$, 
we obtain 
\begin{equation} 
d(\omega^2_1 +\omega^4_3 ) 
+(\omega^2_1 +\omega^4_3 )\wedge da =0 
\label{omega2143} 
\end{equation} 
(using the third and fourth equations in \eqref{domegaij} 
and $\omega^1_2 +\omega^1_4 +\omega^3_2 +\omega^3_4 =0$, 
we can obtain \eqref{omega2143} again).  
From \eqref{omega2143}, we see that there exists a function $b$ 
satisfying 
\begin{equation} 
 \omega^2_1 +\omega^4_3 
=\omega^4_1 +\omega^3_2 
=e^a db. 
\label{beta}
\end{equation} 
From \eqref{dalpha} and \eqref{beta}, 
we obtain 
\begin{equation} 
\omega^3_2 =e^a db -\omega^4_1 , \quad 
\omega^4_3 =e^a db -\omega^2_1 , \quad 
\omega^4_2 =\omega^3_1 -da. 
\label{omega324342} 
\end{equation} 
Applying \eqref{omega324342} to \eqref{domegaij}, we obtain 
\begin{equation} 
\begin{split} 
& d\omega^2_1 +2  \omega^3_1 \wedge  \omega^4_1 
              +e^a db        \wedge  \omega^3_1 
              -da            \wedge  \omega^4_1              =0, \\ 
& d\omega^3_1 +2  \omega^2_1 \wedge  \omega^4_1 
              +e^a db        \wedge (\omega^2_1 -\omega^4_1 )=0, \\ 
& d\omega^4_1 -2  \omega^2_1 \wedge  \omega^3_1 
              -da            \wedge  \omega^2_1 
              +e^a db        \wedge  \omega^3_1              =0. 
\end{split} 
\label{omega213141} 
\end{equation} 
Therefore the connection form $\omega =[\omega^i_j ]$ satisfies 
\eqref{omega213141} for functions $a$, $b$ on $M$. 
Conversely, if $1$-forms $\omega^2_1$, $\omega^3_1$, $\omega^4_1$ 
and functions $a$, $b$ on $M$ satisfy \eqref{omega213141}, 
then $\omega =[\omega^i_j ]$ 
defined by \eqref{Dijn=1} and \eqref{omega324342} 
satisfies $d\omega +\omega \wedge \omega =0$ 
and therefore there exists an ordered pseudo-orthonormal frame field 
$e=(e_1 , e_2 , e_3 , e_4 )$ of $(E, h)$ 
giving the orientation of $E$ and 
satisfying $\nabla eI'_{4, \varepsilon} =eI'_{4, \varepsilon} \omega$. 
Then $e$ defines an $\varepsilon$-nilpotent structure $N$ of $E$ 
such that $e$ is an admissible frame field of $N$. 
Since $\omega^i_j$, $a$, $b$ satisfy \eqref{beta}, 
$N$ satisfies the Walker condition. 
In particular, 
\begin{itemize} 
\item[{\rm (a)}]{$a$ is constant 
if and only if $\nabla N=0$, that is, $\hat{\nabla} \xi_N =0$;} 
\item[{\rm (b)}]{$b$ is constant 
if and only if there exists an $\varepsilon$-nilpotent structure $N'$ 
of $E$ satisfying the Walker condition 
such that $e'$ is an admissible frame field of $N'$.} 
\end{itemize} 
Suppose that $b$ is constant. 
Then \eqref{omega213141} is represented as 
\begin{equation} 
\begin{split} 
& d\omega^2_1 +2\omega^3_1 \wedge \omega^4_1 
              -da          \wedge \omega^4_1 =0, \\ 
& d\omega^3_1 +2\omega^2_1 \wedge \omega^4_1 =0, \\ 
& d\omega^4_1 -2\omega^2_1 \wedge \omega^3_1 
              -da          \wedge \omega^2_1 =0. 
\end{split} 
\label{omega2131412} 
\end{equation} 
If $\omega^3_1 =dc$ for a function $c$ on $M$, 
then 
\begin{equation} 
\omega^2_1 =\dfrac{1}{2} (e^{a  -2c} d\phi 
                         +e^{-a +2c} d\psi ), \quad 
\omega^4_1 =\dfrac{1}{2} (e^{a  -2c} d\phi 
                         -e^{-a +2c} d\psi )
\label{omega2141} 
\end{equation} 
for functions $\phi$, $\psi$ on $M$ 
satisfying $d\phi \wedge d\psi =0$. 
Conversely, 
for functions $c$, $\phi$, $\psi$ on $M$ with $d\phi \wedge d\psi =0$, 
$1$-forms $\omega^2_1$, $\omega^4_1$ defined as in \eqref{omega2141} 
and $\omega^3_1 :=dc$ satisfy \eqref{omega2131412} 
(see Example~\ref{ex:Wnpn=1} 
for the case where $\phi =\psi =0$ and $da\not= 0$). 
\end{ex} 

\begin{ex}\label{ex:dxi0dNnot0}  
Suppose $n\geq 2$. 
Let $M$, $\nabla$ be as in Example~\ref{ex:Wnp}. 
Let $a_1 , \dots , a_n , b_1 , \dots , b_n$ be functions on $M$ 
satisfying 
\begin{equation} 
 \sum^n_{i=1} a_i 
=\sum^n_{i=1} b_i . 
 \label{aibi} 
\end{equation} 
We set 
\begin{equation*} 
A:=\left[ \begin{array}{ccccc} 
           a_1   &  0     & \cdots &  0     \\ 
           0     &  a_2   & \ddots & \vdots \\ 
          \vdots & \ddots & \ddots &  0     \\ 
           0     & \cdots &  0     &  a_n 
            \end{array} 
   \right] , \quad 
B:=\left[ \begin{array}{ccccc} 
           b_1   &  0     & \cdots &  0     \\ 
           0     &  b_2   & \ddots & \vdots \\ 
          \vdots & \ddots & \ddots &  0     \\ 
           0     & \cdots &  0     &  b_n 
            \end{array} 
   \right] 
\end{equation*} 
and 
\begin{equation*} 
 \omega 
=\left[ \begin{array}{cccc} 
          O_n &  O_n & dA   &  O_n \\ 
          O_n &  O_n &  O_n & dB   \\ 
         dA   &  O_n &  O_n &  O_n \\ 
          O_n & dB   &  O_n &  O_n 
          \end{array} 
 \right] . 
\end{equation*} 
Then we have $d\omega +\omega \wedge \omega =O_{4n}$. 
Therefore there exists an ordered pseudo-orthonormal frame field 
$e=(e_1 , \dots , e_{2n}, e_{2n+1} , \dots , e_{4n} )$ of $(E, h)$ 
giving the orientation of $E$ and 
satisfying $\nabla eI'_{4n, \varepsilon} =eI'_{4n, \varepsilon} \omega$. 
Let $N$ be an $\varepsilon$-nilpotent structure of $(E, h)$ 
such that $e$ is an admissible frame field of $N$, 
and $N'$   an $\varepsilon$-nilpotent structure of $(E, h)$ 
such that $e'$ is an admissible frame field of $N'$. 
Then $N$, $N'$ satisfy the Walker condition, by \eqref{NN'W}. 
By \eqref{alpha} and \eqref{aibi}, 
we have $\hat{\nabla} \xi_N    =0$ and 
        $\hat{\nabla} \xi_{N'} =0$. 
If $dA\not= dB$, 
then by \eqref{omegaLambda}, 
$N$, $N'$ are not parallel with respect to $\nabla$. 
\end{ex} 

\begin{ex} 
Let $M$, $\nabla$ be as in Example~\ref{ex:Wnp}. 
Let $F$ be an $n\times n$ symmetric matrix such that each component is 
a function on $M$ and suppose $dF\not= O_n$ and $dF\wedge dF =O_n$. 
We set 
\begin{equation} 
 \omega 
=\left[ \begin{array}{cccc} 
          O_n & O_n & O_n & dF  \\ 
          O_n & O_n & dF  & O_n \\ 
          O_n & dF  & O_n & O_n \\ 
          dF  & O_n & O_n & O_n 
          \end{array} 
 \right] . 
\label{dF} 
\end{equation} 
Then we have $d\omega +\omega \wedge \omega =O_{4n}$. 
Therefore there exists an ordered pseudo-orthonormal frame field 
$e=(e_1 , \dots , e_{2n}, e_{2n+1} , \dots , e_{4n} )$ of $(E, h)$ 
giving the orientation of $E$ and 
satisfying $\nabla eI'_{4n, \varepsilon} =eI'_{4n, \varepsilon} \omega$. 
Let $N$ be an $\varepsilon$-nilpotent structure of $(E, h)$ 
such that $e$ is an admissible frame field of $N$. 
Then $N$ does not satisfy the Walker condition, 
because $\omega$ in \eqref{dF} does not satisfy \eqref{walker}. 
We have 
\begin{equation*} 
(\nabla N)eI'_{4n, \varepsilon}
        =2eI'_{4n, \varepsilon} 
         \left[ \begin{array}{cccc} 
                 O_n &  O_n & dF  &  O_n \\ 
                 O_n &  O_n & O_n & -dF  \\ 
                 dF  &  O_n & O_n &  O_n \\ 
                 O_n & -dF  & O_n &  O_n 
                  \end{array} 
         \right] . 
\end{equation*} 
In particular, if $dF=(1/2)dfI_n$ for a function $f$, 
then $\nabla N=df\otimes J$, 
where $J$ is an $\varepsilon$-paracomplex structure of $(E, h)$ 
defined by \eqref{J}, 
and $\hat{\nabla} \Theta_N =df\otimes \Theta_J$. 
Therefore, if we suppose $n=1$, 
then $\hat{\nabla} \Omega_N =df\otimes \Omega_J$, 
where $\Omega_N$ is the section 
of $U_0\!\left(\bigwedge^2_{\varepsilon}\!E\right)$ 
corresponding to $N$ and given by $(1/\sqrt{2} )\xi_N$, 
and $\Omega_J$ is the section 
of $U_-\!\left(\bigwedge^2_{\varepsilon}\!E\right)$ 
corresponding to $J$. 
\end{ex} 

\section{Special nilpotent structures}\label{sect:sns} 

\setcounter{equation}{0} 

Let $E$ be an oriented vector bundle over $M$ of rank $4n$ 
and $h$ a neutral metric of $E$. 

\begin{defn}\label{defn:K} 
Let $N$ be an $\varepsilon$-nilpotent structure of $(E, h)$. 
Then for a Lie subgroup $K$ of $SO(2n, 2n)$, 
$N$ is called a $K$-\textit{nilpotent structure\/} of $(E, h)$ 
if there exist an open cover $\{ U_{\lambda} \}_{\lambda \in \Lambda}$ of $M$ 
and a family $\{ e_{\lambda} \}_{\lambda \in \Lambda}$ of 
admissible frame fields of $N$ satisfying 
\begin{itemize} 
\item[{\rm (i)}]{each $e_{\lambda}$ is defined on $U_{\lambda}$,} 
\item[{\rm (ii)}]{if $U_{\lambda} \cap U_{\mu} \not= \emptyset$ 
for $\lambda$, $\mu \in \Lambda$, 
then $e_{\mu}     I'_{4n, \varepsilon} 
     =e_{\lambda} I'_{4n, \varepsilon} A_{\lambda \mu}$ 
on $U_{\lambda} \cap U_{\mu}$ 
for a function $A_{\lambda \mu}$ valued in $K$.}  
\end{itemize} 
\end{defn} 

Let $N$ be a $K$-nilpotent structure of $(E, h)$. 
By definition, $N$ gives a principal $K$-bundle $P$ over $M$ 
associated with $E$ such that each $e_{\lambda} I'_{4n, \varepsilon}$ is 
a local section of $P$ 
(see \cite{KN} for principal fiber bundles). 
Let $e$ be an admissible frame field of $N$ 
defined on an open set $U$ of $M$. 
Then $e$ is called a $K$-\textit{admissible frame field\/} 
if $eI'_{4n, \varepsilon}$ gives a local section of $P$, that is, 
if for $\lambda \in \Lambda$ with $U\cap U_{\lambda} \not= \emptyset$, 
there exists a $K$-valued function $A_{\lambda}$ on $U\cap U_{\lambda}$ 
satisfying $e           I'_{4n, \varepsilon} 
           =e_{\lambda} I'_{4n, \varepsilon} A_{\lambda}$. 
Based on Definition~\ref{defn:K}, we obtain 

\begin{pro}\label{pro:G} 
Any $\varepsilon$-nilpotent structure of $(E, h)$ is 
a $G$-nilpotent structure of $(E, h)$. 
\end{pro} 

Let $\nabla$ be an $h$-connection of $E$. 

\begin{defn}\label{defn:K2} 
Let $N$ be an $\varepsilon$-nilpotent structure of $(E, h)$. 
Then for a Lie subgroup $K$ of $SO(2n, 2n)$, 
$N$ is called a $K$-\textit{nilpotent structure\/} of $(E, h, \nabla )$ 
if there exist an open cover $\{ U_{\lambda} \}_{\lambda \in \Lambda}$ of $M$ 
and a family $\{ e_{\lambda} \}_{\lambda \in \Lambda}$ of 
admissible frame fields of $N$ 
satisfying (i), (ii) in Definition~\ref{defn:K} and 
\begin{itemize} 
\item[{\rm (iii)}]{the connection form $\omega_{\lambda}$ of $\nabla$ 
with respect to $e_{\lambda} I'_{4n, \varepsilon}$ is valued in 
the Lie algebra of $K$.}  
\end{itemize} 
\end{defn} 

Let $N$ be a $K$-nilpotent structure of $(E, h, \nabla )$. 
Then $\nabla$ gives a connection in the principal $K$-bundle $P$. 
Let $e$ be a $K$-admissible frame field of $N$. 
Then the connection form $\omega$ of $\nabla$ 
with respect to $eI'_{4n, \varepsilon}$ 
is given by $\omega =A^{-1}_{\lambda} \omega_{\lambda} A_{\lambda} 
                    +A^{-1}_{\lambda}                 dA_{\lambda}$. 
Therefore $\omega$ is valued in the Lie algebra of $K$. 

Suppose that $\nabla$ is flat, that is, 
the curvature tensor of $\nabla$ vanishes. 
Then there exists an ordered pseudo-orthonormal local frame field $e$ 
of $(E, h)$ on a neighborhood $U$ of each point of $M$ 
giving the orientation of $E$ 
such that each local section of $E$ which appears in $e$ is parallel 
with respect to $\nabla$. 
Then the connection form of $\nabla$ with respect to $e$ vanishes. 
Therefore for $\lambda \in \Lambda$ 
with $U\cap U_{\lambda} \not= \emptyset$, 
if we define an $SO(2n, 2n)$-valued function $\Tilde{A}_{\lambda}$ 
on $U\cap U_{\lambda}$ 
by $e           I'_{4n, \varepsilon} 
   =e_{\lambda} I'_{4n, \varepsilon} \Tilde{A}_{\lambda}$, 
then $\Tilde{A}_{\lambda}$ satisfies 
$d\Tilde{A}_{\lambda} =-\omega_{\lambda} \Tilde{A}_{\lambda}$. 
This relation means that $d\Tilde{A}_{\lambda}$ is valued in 
the distribution $\mathscr{D}$ on $SO(2n, 2n)$ 
obtained by the right translations of the Lie algebra of $K$. 
If we suppose that $U\cap U_{\lambda}$ is connected, 
then there exists an integral manifold of $\mathscr{D}$ 
containing $\Tilde{A}_{\lambda}$. 
We can choose $e$ so that $\Tilde{A}_{\lambda}$ is valued in $K$ 
at a point. 
Since $K$ is an integral manifold of $\mathscr{D}$, 
$\Tilde{A}_{\lambda}$ is valued in $K$ over its domain. 
If we suppose $K\subset G$, 
then $e$ is a $K$-admissible frame field of $N$. 
Hence we obtain 

\begin{pro}\label{pro:Kflat} 
Suppose that $\nabla$ is flat. 
Let $K$ be a Lie subgroup of $G$. 
Then for a $K$-nilpotent structure $N$ of $(E, h, \nabla )$, 
there exists a $K$-admissible frame field $e$ of $N$ 
on a neighborhood of each point of $M$ 
such that each local section of $E$ which appears in $e$ is parallel 
with respect to $\nabla$. 
\end{pro} 

\begin{rem} 
Let $\eta_1 , \dots , \eta_n , \eta_{n+1} , \dots , \eta_{2n}$ 
be $2n$ vectors of $E^{4n}_{2n}$ defined by 
\begin{equation*} 
\Lambda_n =[\eta_{n+1} \ \dots \  \eta_{2n} \ 
           -\eta_1     \ \dots \ -\eta_n    \ 
            \eta_{n+1} \ \dots \  \eta_{2n} \ 
            \eta_1     \ \dots \  \eta_n ]. 
\end{equation*} 
Let $W$ be a $2n$-dimensional subspace of $E^{4n}_{2n}$ 
spanned by $\eta_1 , \dots , \eta_n , \eta_{n+1} , \dots , \eta_{2n}$. 
Then $W$ is a light-like subspace of $E^{4n}_{2n}$. 
Let $SO(2n, 2n)_W$ be a Lie subgroup of $SO(2n, 2n)$ 
which consists of elements preserving $W$. 
Then we have $G\subset SO(2n, 2n)_W$ 
and $G\not= SO(2n, 2n)_W$ (\cite{ando8}). 
An $SO(2n, 2n)_W$-nilpotent structure of $(E, h, \nabla )$ 
is just an $\varepsilon$-nilpotent structure of $(E, h)$ 
satisfying the Walker condition. 
If $K=SO(2n, 2n)_W$, 
then the conclusion of Proposition~\ref{pro:Kflat} is not necessarily valid, 
because of Example~\ref{ex:Wnp} $\sim$ Example~\ref{ex:dxi0dNnot0}, 
while for an almost $\varepsilon$-nilpotent structure $N$ of $E^{4n}_{2n}$ 
satisfying the Walker condition 
with respect to the Levi-Civita connection, 
$\pi_N$ gives a constant $2n$-dimensional light-like subspace 
of $E^{4n}_{2n}$. 
\end{rem} 

We will prove 

\begin{pro}\label{pro:G2} 
A $G$-nilpotent structure of $(E, h, \nabla )$ is just 
an $\varepsilon$-nilpotent structure of $(E, h)$ 
parallel with respect to $\nabla$. 
\end{pro} 

\vspace{3mm} 

\par\noindent 
\textit{Proof} \ 
Let $N$ be a $G$-nilpotent structure of $(E, h, \nabla )$. 
Then the connection form $\omega$ of $\nabla$ with respect to 
$eI'_{4n, \varepsilon}$ for any admissible frame field $e$ of $N$ 
is valued in the Lie algebra of $G$. 
Therefore we have $\omega \Lambda_n =\Lambda_n \omega$. 
By this, we obtain $\nabla (Ne)=N\nabla e$. 
Therefore $N$ is parallel with respect to $\nabla$. 
Suppose that $N$ is an $\varepsilon$-nilpotent structure of $(E, h)$ 
parallel with respect to $\nabla$. 
Then we have $\omega \Lambda_n =\Lambda_n \omega$ 
for the connection form $\omega$ of $\nabla$ with respect to 
$eI'_{4n, \varepsilon}$ for any admissible frame field $e$ of $N$. 
This means that $\omega$ is valued in the Lie algebra of $G$. 
Therefore we can find an open cover $\{ U_{\lambda} \}_{\lambda \in \Lambda}$ 
of $M$ and a family $\{ e_{\lambda} \}_{\lambda \in \Lambda}$ of 
admissible frame fields of $N$ 
satisfying (i), (ii) in Definition~\ref{defn:K} and 
(iii)                in Definition~\ref{defn:K2} for $K=G$. 
Therefore $N$ is a $G$-nilpotent structure of $(E, h, \nabla )$. 
\hfill 
$\square$ 

\vspace{3mm} 

Let $H$ be a Lie subgroup of $SO(2n, 2n)$ defined by 
\begin{equation} 
H:=\left\{ \left[ 
           \begin{array}{cccc} 
            A_{11} & -A_{21} &  A_{31} & A_{41} \\ 
            A_{21} &  A_{11} & -A_{41} & A_{31} \\ 
            A_{31} & -A_{41} &  A_{11} & A_{21} \\ 
            A_{41} &  A_{31} & -A_{21} & A_{11} 
             \end{array} 
           \right] \in SO(2n, 2n) 
   \right\} . 
\label{H} 
\end{equation} 
Then $H$ is a Lie subgroup of $G$. 
If $N$ is an $H$-nilpotent structure of $(E, h, \nabla )$, 
then for an $H$-admissible frame field $e$ of $N$, 
the connection form $\omega$ of $\nabla$ 
with respect to $eI'_{4n, \varepsilon}$ is represented as 
\begin{equation}
 \omega 
=\left[ \begin{array}{cccc} 
         D_{11} & -D_{21} &  D_{31} & D_{41} \\ 
         D_{21} &  D_{11} & -D_{41} & D_{31} \\ 
         D_{31} & -D_{41} &  D_{11} & D_{21} \\ 
         D_{41} &  D_{31} & -D_{21} & D_{11} 
          \end{array} 
 \right] \qquad 
 \left( 
 \begin{array}{l} 
  {}^t D_{11} =-D_{11} , \\ 
  {}^t D_{21} = D_{21} , \\ 
  {}^t D_{31} = D_{31} , \\ 
  {}^t D_{41} = D_{41} 
   \end{array} 
 \right) . 
\label{somega} 
\end{equation} 
The main objects of study in this section are $H$-nilpotent structures 
of $(E, h, \nabla )$. 

Let $I$ be a complex structure of $E$ 
and $J_1$, $J_2$ paracomplex structures of $E$. 
Suppose that $I$, $J_1$, $J_2$ give 
a \textit{paraquaternionic structure\/} $E_{I, J_1 , J_2}$, 
that is, 
suppose that $I$, $J_1$, $J_2$ satisfy 
\begin{itemize} 
\item[{\rm (i)}]{$IJ_1 =J_2$} 
\end{itemize} 
and span a subbundle $E_{I, J_1 , J_2}$ of ${\rm End}\,E$ of rank $3$. 
In addition, suppose that $h$ is \textit{adapted\/} to $E_{I, J_1 , J_2}$, 
that is, suppose 
\begin{itemize} 
\item[{\rm (ii)}]{$I$ preserves $h$ and $J_1$, $J_2$ reverse $h$, that is, 
$h$ is Hermitian with respect to $I$ and 
   paraHermitian with respect to $J_1$, $J_2$.} 
\end{itemize} 
See \cite{BDM2}, \cite{GMV} for paraquaternionic structures 
(notice the sign in (i)). 
We say that $h$, $\nabla$, $I$, $J_1$, $J_2$ form 
a \textit{neutral hyperK\"{a}hler structure\/} of $E$ 
if $h$, $\nabla$, $I$, $J_1$, $J_2$ satisfy (i), (ii) and 
\begin{itemize} 
\item[{\rm (iii)}]{$I$, $J_1$, $J_2$ are parallel with respect to $\nabla$.} 
\end{itemize} 
See \cite{DGMY}, \cite{kamada} for neutral hyperK\"{a}hler $4$-manifolds. 

We will prove 

\begin{thm}\label{thm:nh} 
The following hold\/$:$ 
\begin{itemize} 
\item[{\rm (a)}]{the section of $E_{I, J_1 , J_2}$ given by 
$N_{r, \theta} :=r(I-(\sin \theta )J_1 +(\cos \theta )J_2 )$ 
for $r\in \mbox{\boldmath{$R$}} \setminus \{ 0\}$ 
and $\theta \in [0, 2\pi )$ is an $H$-nilpotent structure of $(E, h);$} 
\item[{\rm (b)}]{if $h$, $\nabla$, $I$, $J_1$, $J_2$ form 
a neutral hyperK\"{a}hler structure of $E$, 
then $N_{r, \theta}$ is an $H$-nilpotent structure of $(E, h, \nabla )$.} 
\end{itemize} 
\end{thm} 

\vspace{3mm} 

\par\noindent 
\textit{Proof} \ 
Let $\{ U_{\lambda} \}_{\lambda \in \Lambda}$ be an open cover of $M$ 
such that on each $U_{\lambda}$, 
there exists an ordered pseudo-orthonormal local frame field 
$e_{\lambda} =(e_{\lambda ,    1} , \dots , e_{\lambda , 2n}, 
               e_{\lambda , 2n+1} , \dots , e_{\lambda , 4n} )$ of $E$ 
satisfying 
\begin{equation} 
I   e_{\lambda , i} =e_{\lambda ,  n+i} , \quad 
J_1 e_{\lambda , i} =e_{\lambda , 2n+i} , \quad 
J_2 e_{\lambda , i} =e_{\lambda , 3n+i}   \quad 
(i=1, \dots , n). 
\label{IJ1J2e} 
\end{equation} 
If $U_{\lambda} \cap U_{\mu} \not= \emptyset$ 
for $\lambda$, $\mu \in \Lambda$ and 
if we represent $e_{\mu}$ as $e_{\mu} =e_{\lambda} A_{\lambda \mu}$ 
on $U_{\lambda} \cap U_{\mu}$, 
then $A_{\lambda \mu}$ is valued in $H$. 
Hence $I$, $J_1$, $J_2$ define a principal $H$-bundle $P$ over $M$ 
associated with $E$ so that each $e_{\lambda}$ is a local section of $P$. 
For each $\lambda \in \Lambda$, 
let $N_{\lambda} $ be a section of 
the restriction of ${\rm End}\,E$ on $U_{\lambda}$ 
defined by $N_{\lambda} e_{\lambda} =e_{\lambda} \Lambda_n$. 
Then noticing that $H$ is contained in $G$, 
we can define a section $N$ of ${\rm End}\,E$ by $N=N_{\lambda}$ 
on each $U_{\lambda}$, which is an $\varepsilon$-nilpotent structure 
of $(E, h)$ for $\varepsilon =+$ or $-$ 
according to whether $e_{\lambda}$ gives the orientation of $E$ or not. 
In addition, $N$ is an $H$-nilpotent structure of $(E, h)$ 
and for each $\lambda \in \Lambda$, 
$e_{\lambda} I'_{4n, \varepsilon}$ is an $H$-admissible frame field of $N$. 
By \eqref{IJ1J2e}, $N$ coincides with $I+J_2 =N_{1, 0}$. 
For $\theta \in [0, 2\pi )$, 
we set 
\begin{equation*} 
J'_1 := (\cos \theta )J_1 +(\sin \theta )J_2 , \quad 
J'_2 :=-(\sin \theta )J_1 +(\cos \theta )J_2 . 
\end{equation*} 
Then $J'_1$, $J'_2$ are paracomplex structures of $E$ reversing $h$, and 
$I$, $J'_1$, $J'_2$ satisfy $IJ'_1 =J'_2$ and span $E_{I, J_1 , J_2}$. 
Therefore referring to the above discussion, 
we obtain an $H$-nilpotent structure of $(E, h)$, 
which is given by 
\begin{equation*} 
I+J'_2 =I-(\sin \theta )J_1 +(\cos \theta )J_2 
       =N_{1, \theta} . 
\end{equation*} 
For $t\not= 0$ and $\delta \in \{ 1, -1\}$, 
we set 
\begin{equation*} 
I'    :=\delta (({\rm cosh}\,t)I +({\rm sinh}\,t)J'_2 ), \quad 
J''_2 :=\delta (({\rm sinh}\,t)I +({\rm cosh}\,t)J'_2 ). 
\end{equation*} 
Then $I'$    is a     complex structure of $E$ preserving $h$ 
and  $J''_2$ is a paracomplex structure of $E$ reversing $h$. 
In addition, $I'$, $J'_1$ and $J''_2$ satisfy $I'J'_1 =J''_2$ 
and span $E_{I, J_1 , J_2}$. 
Therefore referring to the above discussion, 
we obtain an $H$-nilpotent structure of $(E, h)$ 
given by 
\begin{equation*} 
I' +J''_2 =\delta e^t (I-(\sin \theta )J_1 +(\cos \theta )J_2 ) 
          =N_{\delta e^t , \theta} . 
\end{equation*} 
Hence we obtain (a) in Theorem~\ref{thm:nh}. 
We set 
\begin{equation*} 
\begin{split} 
  (\nabla e_{\lambda , 1} \ \dots \ \nabla e_{\lambda , n} )= 
&  (e_{\lambda ,    1} \ \dots \ e_{\lambda ,  n} )D_{11} 
  +(e_{\lambda ,  n+1} \ \dots \ e_{\lambda , 2n} )D_{21} \\ 
& +(e_{\lambda , 2n+1} \ \dots \ e_{\lambda , 3n} )D_{31} 
  +(e_{\lambda , 3n+1} \ \dots \ e_{\lambda , 4n} )D_{41} . 
\end{split} 
\end{equation*} 
Suppose that $h$, $\nabla$, $I$, $J_1$, $J_2$ form 
a neutral hyperK\"{a}hler structure of $E$. 
Then by $\nabla I=0$, we obtain 
\begin{equation*} 
\begin{split} 
  (\nabla e_{\lambda , n+1} \ \dots \ \nabla e_{\lambda , 2n} )=
& -(e_{\lambda ,    1} \ \dots \ e_{\lambda ,  n} )D_{21} 
  +(e_{\lambda ,  n+1} \ \dots \ e_{\lambda , 2n} )D_{11} \\ 
& -(e_{\lambda , 2n+1} \ \dots \ e_{\lambda , 3n} )D_{41} 
  +(e_{\lambda , 3n+1} \ \dots \ e_{\lambda , 4n} )D_{31} . 
\end{split} 
\end{equation*} 
Similarly, by $\nabla J_1 =0$, we obtain 
\begin{equation*} 
\begin{split} 
  (\nabla e_{\lambda , 2n+1} \ \dots \ \nabla e_{\lambda , 3n} )= 
&  (e_{\lambda ,    1} \ \dots \ e_{\lambda ,  n} )D_{31} 
  -(e_{\lambda ,  n+1} \ \dots \ e_{\lambda , 2n} )D_{41} \\ 
& +(e_{\lambda , 2n+1} \ \dots \ e_{\lambda , 3n} )D_{11} 
  -(e_{\lambda , 3n+1} \ \dots \ e_{\lambda , 4n} )D_{21} , 
\end{split} 
\end{equation*} 
and by $\nabla J_2 =0$, we obtain 
\begin{equation*} 
\begin{split} 
  (\nabla e_{\lambda , 3n+1} \ \dots \ \nabla e_{\lambda , 4n} )= 
&  (e_{\lambda ,    1} \ \dots \ e_{\lambda ,  n} )D_{41} 
  +(e_{\lambda ,  n+1} \ \dots \ e_{\lambda , 2n} )D_{31} \\ 
& +(e_{\lambda , 2n+1} \ \dots \ e_{\lambda , 3n} )D_{21} 
  +(e_{\lambda , 3n+1} \ \dots \ e_{\lambda , 4n} )D_{11} . 
\end{split} 
\end{equation*} 
Therefore the connection form $\omega_{\lambda}$ of $\nabla$ 
with respect to $e_{\lambda}$ is given as in \eqref{somega} 
and therefore the connection form $\omega_{\lambda}$ is valued in 
the Lie algebra of $H$. 
Hence $\nabla$ gives a connection in $P$ 
and $N=N_{1, 0}$ is an $H$-nilpotent structure of $(E, h, \nabla )$. 
Similarly, 
for $r\in \mbox{\boldmath{$R$}} \setminus \{ 0\}$ 
and $\theta \in [0, 2\pi )$, 
$N_{r, \theta}$ is an $H$-nilpotent structure of $(E, h, \nabla )$. 
Hence we obtain (b) in Theorem~\ref{thm:nh}. 
\hfill 
$\square$ 

\vspace{3mm} 

\begin{rem} 
Whether $e_{\lambda}$ gives the orientation of $E$ does not depend 
on $\lambda \in \Lambda$, 
and $N_{r, \theta}$ is a $+$ or $-$-nilpotent structure 
according to whether $e_{\lambda}$ gives the orientation or not, 
as was already seen in the above proof. 
\end{rem} 

Let $N$ be an $\varepsilon$-nilpotent structure of $(E, h)$ 
and suppose that $N$ is an $H$-nilpotent structure of $(E, h)$. 
Then the light-like subbundle $\pi_N$ of rank $2n$ is determined by $N$ 
and locally spanned by $\xi_1 , \dots , \xi_{2n}$ as in \eqref{xie} 
for an $H$-admissible frame field $e$ of $N$. 
In addition, we will prove 

\begin{thm}\label{thm:pi'N'} 
Let $N$ be as above. 
\begin{itemize} 
\item[{\rm (a)}]{There exists a unique light-like subbundle $\pi^{\times}_N$ 
of $(E, h)$ of rank $2n$ which is locally spanned by 
\begin{equation} 
  \begin{array}{lcl} 
   \xi'_1      :=e_1     +            e_{2n+1} , & \ &  
   \xi'_i      :=e_i     +            e_{2n+i} , \\ 
   \xi'_{n+1}  :=e_{n+1} -\varepsilon e_{3n+1} , & \ &  
   \xi'_{n+i}  :=e_{n+i} -            e_{3n+i}  
    \end{array} \ 
(i=2, \dots , n) 
\label{xie'} 
\end{equation} 
for any $H$-admissible frame field $e$ of $N$. 
In particular, $E=\pi_N \oplus \pi^{\times}_N$.} 
\item[{\rm (b)}]{There exists 
a unique $\varepsilon$-nilpotent and $H$-nilpotent structure $N^{\times}$ 
of $(E, h)$ such that each $H$-admissible frame field of $N^{\times}$ is 
given by $e' =(e_1 , \dots , e_{2n}, -e_{2n+1} , \dots , -e_{4n} )$ 
for an $H$-admissible frame field $e$ of $N$. 
In particular, $\pi_{N^{\times}} =\pi^{\times}_N$.} 
\item[{\rm (c)}]{The condition that $N$ is an $H$-nilpotent structure 
of $(E, h, \nabla )$ is equivalent to 
the condition that $N^{\times}$ is an $H$-nilpotent structure 
of $(E, h, \nabla )$.} 
\end{itemize} 
\end{thm} 

\vspace{3mm} 

\par\noindent 
\textit{Proof} \ 
For an $H$-admissible frame field $e$ of $N$, 
the subbundle of the restriction of $E$ on the domain of $e$ 
spanned by $\xi'_i$, $\xi'_{n+i}$ ($i=1, \dots , n$) as in \eqref{xie'} 
does not depend on the choice of $e$, 
because of the definition of $H$ in \eqref{H}. 
Therefore there exists a unique light-like subbundle $\pi^{\times}_N$ 
of $(E, h)$ of rank $2n$ satisfying (a) in Theorem~\ref{thm:pi'N'}. 
Hence we obtain (a) in Theorem~\ref{thm:pi'N'}. 
For $N$, we can find an open cover $\{ U_{\lambda} \}_{\lambda \in \Lambda}$ 
of $M$ and a family $\{ e_{\lambda} \}_{\lambda \in \Lambda}$ of 
$H$-admissible frame fields of $N$ satisfying (i), (ii) 
in Definition~\ref{defn:K} for $K=H$. 
We set 
\begin{equation*} 
A'_{\lambda \mu} :=I_{2n, 2n} A_{\lambda \mu} I_{2n, 2n} , \quad 
I_{2n, 2n} 
:=\left[ \begin{array}{cccc} 
          I_n & O_n &  O_n &  O_n \\ 
          O_n & I_n &  O_n &  O_n \\ 
          O_n & O_n & -I_n &  O_n \\ 
          O_n & O_n &  O_n & -I_n 
           \end{array} 
  \right] 
\end{equation*} 
for $\lambda$, $\mu \in \Lambda$ 
satisfying $U_{\lambda} \cap U_{\mu} \not= \emptyset$. 
Then $A'_{\lambda \mu}$ satisfies 
$e'_{\mu}     I'_{4n, \varepsilon} 
=e'_{\lambda} I'_{4n, \varepsilon}A'_{\lambda \mu}$ 
on $U_{\lambda} \cap U_{\mu}$, and 
since $A_{\lambda \mu}$ is valued in $H$, 
$A'_{\lambda \mu}$ is also valued in $H$. 
Therefore there exists an $\varepsilon$-nilpotent and $H$-nilpotent 
structure $N^{\times}$ of $(E, h)$ 
such that $e'_{\lambda}$ ($\lambda \in \Lambda$) are 
$H$-admissible frame fields of $N^{\times}$ 
and then $N^{\times}$ is uniquely determined 
by $\{ e'_{\lambda} \}_{\lambda \in \Lambda}$. 
Let $\Tilde{e}$ be an $H$-admissible frame field of $N^{\times}$ 
defined on an open set $U$ of $M$. 
Then for $\lambda \in \Lambda$ with $U\cap U_{\lambda} \not= \emptyset$, 
there exists an $H$-valued function $A'_{\lambda}$ on $U\cap U_{\lambda}$ 
satisfying 
$\Tilde{e}    I'_{4n, \varepsilon} 
=e'_{\lambda} I'_{4n, \varepsilon} A'_{\lambda}$. 
This relation is equivalent to 
$\Tilde{e}'  I'_{4n, \varepsilon} 
=e_{\lambda} I'_{4n, \varepsilon} A_{\lambda}$ 
with $A_{\lambda} :=I_{2n, 2n} A'_{\lambda} I_{2n, 2n}$. 
Since $A_{\lambda}$ is $H$-valued, 
$\Tilde{e}'$ is an $H$-admissible frame field of $N$. 
Hence we obtain (b) in Theorem~\ref{thm:pi'N'}. 
Suppose that $N$ is an $H$-nilpotent structure of $(E, h, \nabla )$. 
Then we can suppose that the connection form $\omega_{\lambda}$ 
of $\nabla$ with respect to $e_{\lambda} I'_{4n, \varepsilon}$ is valued in 
the Lie algebra of $H$. 
Therefore the connection form $\omega'_{\lambda}$ of $\nabla$ 
with respect to $e'_{\lambda} I'_{4n, \varepsilon}$ is valued 
in the Lie algebra of $H$. 
This means that $N^{\times}$ is an $H$-nilpotent structure 
of $(E, h, \nabla )$. 
Similarly, we can prove the converse. 
Hence we obtain (c) in Theorem~\ref{thm:pi'N'}. 
\hfill 
$\square$ 

\vspace{3mm} 

Let $N$ be an $H$-nilpotent structure of $(E, h)$. 
Then we call $N^{\times}$ as in Theorem~\ref{thm:pi'N'} 
the \textit{dual\/} $H$-\textit{nilpotent structure\/} of $N$. 
We see that $N$ is the dual $H$-nilpotent structure of $N^{\times}$: 
$(N^{\times} )^{\times} =N$. 

For an $H$-nilpotent structure $N$ of $(E, h)$, 
we set 
\begin{equation} 
I   :=\dfrac{1}{2} (N+N^{\times} ), \quad 
J_2 :=\dfrac{1}{2} (N-N^{\times} ), \quad 
J_1 :=-IJ_2 . 
\label{IJ2J1} 
\end{equation} 
Then $I$, $J_1$, $J_2$ give 
a paraquaternionic structure $E_{I, J_1 , J_2}$ 
such that $h$ is adapted to $E_{I, J_1 , J_2}$. 
In addition, 
if $N$ is an $H$-nilpotent structure of $(E, h, \nabla )$, 
then $h$, $\nabla$, $I$, $J_1$, $J_2$ form 
a neutral hyperK\"{a}hler structure of $E$. 
Hence we obtain 

\begin{cor}\label{cor:pi'N'} 
An $H$-nilpotent structure $N$ of $(E, h)$ defines 
a paraquaternionic structure $E_{I, J_1 , J_2}$ 
such that $h$ is adapted to $E_{I, J_1 , J_2}$ by \eqref{IJ2J1}. 
In addition, 
if $N$ is an $H$-nilpotent structure of $(E, h, \nabla )$, 
then $h$, $\nabla$, $I$, $J_1$, $J_2$ form 
a neutral hyperK\"{a}hler structure of $E$. 
\end{cor} 

\begin{rem} 
Let $E$, $h$ be as in the beginning of Section~\ref{sect:ans}. 
Let $I$ be a section of ${\rm End}\,E$. 
We say that $I$ is an $\varepsilon$-\textit{complex structure\/} 
of $(E, h)$ if $I$ satisfies 
\begin{itemize} 
\item[{\rm (i)}]{$I$ is a complex structure of $E$,}
\item[{\rm (ii)}]{$I$ is $h$-preserving, that is, $I^* h=h$,} 
\item[{\rm (iii)}]{on a neighborhood of each point of $M$, 
there exists an ordered pseudo-orthonormal local frame field 
$e=(e_1 ,\dots e_{2n} , e_{2n+1} , \dots , e_{4n} )$ of $E$ 
giving the orientation and satisfying 
\begin{equation*} 
IeI'_{4n, \varepsilon} 
=eI'_{4n, \varepsilon} \Lambda_{n, +} , \quad 
  \Lambda_{n, +} 
:=\left[ \begin{array}{cccc} 
          O_n & -I_n & O_n &  O_n \\ 
          I_n &  O_n & O_n &  O_n \\ 
          O_n &  O_n & O_n & -I_n \\ 
          O_n &  O_n & I_n &  O_n 
           \end{array} 
  \right] . 
\label{I} 
\end{equation*}} 
\end{itemize}  
Let $I$ be an $\varepsilon$-complex structure of $(E, h)$. 
Then such a frame field as $e$ is called 
an \textit{admissible frame field\/} of $I$. 
We see that $I$ gives a section $I^*$ of ${\rm End}\,E^*$ 
by $I^* \phi^* =\phi^* \circ I$ 
for a local section $\phi^*$ of $E^*$. 
Then we have 
$I^* e^* I'_{4n, \varepsilon} 
    =e^* I'_{4n, \varepsilon} {}^t\!\Lambda_{n, +}$. 
Therefore, 
if we set $\Theta_I (\phi^* , \psi^* ):=h^* (\phi^* , I^* \psi^* )$, 
then $\Theta_I$ is a section of $\bigwedge^2\!E$ and 
locally represented as 
\begin{equation} 
 \Theta_I 
=\sum^n_{i=1} (e_i      \wedge e_{n+i} 
                      -(\varepsilon 1)^{\delta_{i1}} 
               e_{2n+i} \wedge e_{3n+i} ). 
\label{ThetaI} 
\end{equation} 
For an $h$-connection $\nabla$ of $E$, 
$\nabla I=0$ is equivalent to $\hat{\nabla} \Theta_I =0$. 
Let $N$ be an $\varepsilon$-nilpotent and $H$-nilpotent structure 
of $(E, h)$ and $e$ an $H$-admissible frame field of $N$. 
Let $I$, $J_2$ be as in \eqref{IJ2J1}. 
Then $I$ is an $\varepsilon$-complex structure of $(E, h)$ 
such that $e$ is an admissible frame field of $I$ 
and $J_2$ is an $\varepsilon$-paracomplex structure of $(E, h)$ 
in the sense of Remark~\ref{rem:J} such that 
\begin{equation*} 
(e_{n+1}  , \dots , e_{2n} , -e_1      , \dots , -e_n , 
 e_{2n+1} , \dots , e_{3n} ,  e_{3n+1} , \dots , e_{4n} ) 
\end{equation*} 
is an admissible frame field of $J_2$. 
In addition, $\Theta_N =\Theta_I +\Theta_{J_2}$. 
Let $J_1$ be as in \eqref{IJ2J1}. 
Then $J_1$ is an $\varepsilon$-paracomplex structure of $(E, h)$ 
in the sense of Remark~\ref{rem:J} such that 
$e$ is an admissible frame field of $J_1$. 
\end{rem} 

\begin{rem} 
Let $I$ be a complex structure of $E$ 
and $J_1$, $J_2$ paracomplex structures of $E$. 
Suppose that $I$, $J_1$, $J_2$ give 
a paraquaternionic structure $E_{I, J_1 , J_2}$ 
such that $h$ is adapted to $E_{I, J_1 , J_2}$. 
Then $N   :=I+J_2$ is an $H$-nilpotent structure of $(E, h)$ 
and  $N^{\times} :=I-J_2$ is its dual $H$-nilpotent structure. 
\end{rem} 

\begin{rem} 
If $N$ is an $H$-nilpotent structure of $(E, h, \nabla )$, 
then we see by \eqref{somega} that 
the connection form of an $H$-admissible frame field $e$ of $N$ 
satisfies \eqref{NN'W}. 
\end{rem} 

\begin{rem} 
Let $N$ be an $\varepsilon$-nilpotent structure of $(E, h)$. 
Let $e$, $\Tilde{e}$ be admissible frame fields of $N$ 
and $A$ a $G$-valued function on the intersection of the domains 
of $e$, $\Tilde{e}$ given by $\Tilde{e} I'_{4n, \varepsilon} 
                             =       e  I'_{4n, \varepsilon} A$. 
Then the subbundle spanned by $\xi_i$, $\xi_{n+i}$ ($i=1, \dots , n$) 
coincides with 
the subbundle spanned by $\Tilde{\xi}_i$, $\Tilde{\xi}_{n+i}$ 
($i=1, \dots , n$). 
In addition, 
if the subbundle spanned by $\xi'_i$, $\xi'_{n+i}$ ($i=1, \dots , n$) 
coincides with 
the subbundle spanned by $\Tilde{\xi}'_i$, $\Tilde{\xi}'_{n+i}$ 
($i=1, \dots , n$), 
then $A$ is valued in $H$. 
\end{rem}

\vspace{4mm} 

\par\noindent 
\footnotesize{Faculty of Advanced Science and Technology, 
              Kumamoto University \\ 
              2--39--1 Kurokami, Chuo-ku, Kumamoto 860--8555 Japan} 

\par\noindent  
\footnotesize{E-mail address: andonaoya@kumamoto-u.ac.jp} 

\end{document}